\documentclass[a4paper,11pt]{article}
\usepackage{amsfonts}
\usepackage{amstext}
\usepackage{amsthm}
\usepackage{amsmath}
\usepackage{amssymb}
\usepackage{latexsym}
\usepackage[latin1]{inputenc}
\newcommand{\R}{\Bbb{R}}
\newcommand{\N}{\Bbb{N}}

\newtheorem{teor}{Theorem}[section]
\newtheorem{propo}{Proposition}[section]

\newtheorem{lema}{Lemma}[section]

\newtheorem{rem}{Remark}[section]

\newcommand{\n}{\noindent}

\newcommand {\fim}{\rule{0.5em}{0.5em}}

\begin{document}

\title{A priori bounds and positive solutions for non-variational fractional elliptic systems
}

\author{\textbf{Edir Junior Ferreira Leite, Marcos Montenegro \footnote{\textit{E-mail addresses}:
edirjrleite@msn.com (E.J.F. Leite), montene@mat.ufmg.br (M.
Montenegro)}}\\ {\small\it Departamento de Matem\'{a}tica,
Universidade Federal de Minas Gerais,}\\ {\small\it Caixa Postal
702, 30123-970, Belo Horizonte, MG, Brazil}}

\date{}{

\maketitle

\markboth{abstract}{abstract}
\addcontentsline{toc}{chapter}{abstract}

\hrule \vspace{0,2cm}

\n {\bf Abstract}

In this paper we study strongly coupled elliptic systems in non-variational form involving fractional Laplace operators. We prove Liouville type theorems and, by mean of the blow-up method, we establish a priori bounds of positive solutions for subcritical and superlinear nonlinearities in a coupled sense. By using those latter, we then derive the existence of positive solutions through topological methods.

\vspace{0.5cm}
\hrule\vspace{0.2cm}

\section{Introduction and main results}

The present paper deals with a priori bounds and existence of positive solutions for elliptic systems of the form

\begin{equation} \label{1}
\left\{
\begin{array}{llll}
(-\Delta)^{s}u = v^p & {\rm in} \ \ \Omega\\
(-\Delta)^{t}v = u^q & {\rm in} \ \ \Omega\\
u= v=0 & {\rm in} \ \ \R^n\setminus\Omega
\end{array}
\right.
\end{equation}
where $\Omega$ is a smooth bounded open subset of $\R^{n}$, $n\geq 2$, $0 <s,t < 1$, $p,q > 0$ and the fractional Laplace operator $(-\Delta)^{s}$ is defined as

\begin{equation} \label{2}
(-\Delta)^{s}u(x) = C(n,s)\, P.V.\int\limits_{\R^{n}}\frac{u(x)-u(y)}{\vert x-y\vert^{n+2s}}\; dy\, ,
\end{equation}
or equivalently,

\[
(-\Delta)^{s}u(x)=-\frac{1}{2}C(n,s)\int\limits_{\R^{n}}\frac{u(x+y)+u(x-y)-2u(x)}{|y|^{n+2s}}\; dy
\]

\n for all $x \in \R^{n}$, where P.V. denotes the principal value of integrals and

\[
C(n,s) = \left(\int\limits_{\R^{n}}\frac{1-\cos(\zeta_{1})}{\vert\zeta\vert^{n+2s}}\; d\zeta\right)^{-1}
\]

\n with $\zeta = (\zeta_1, \ldots, \zeta_n) \in \R^n$.

The operator $(-\Delta)^{s}$ satisfies

\[
\lim_{s\longrightarrow 1^{-}}(-\Delta)^{s} u = -\Delta u
\]

\n pointwise in $\R^n$ for all $u \in C^{\infty}_{0}(\R^{n})$, so it interpolates the Laplace operator in $\R^n$.

A closely related to (but different from) the restricted fractional Laplace operator $(-\Delta)^{s}$ is the spectral fractional Laplace operator $\mathcal{A}^{s}$ with zero Dirichlet boundary values on $\partial \Omega$. Its definition can be given in terms of the Dirichlet spectra of the Laplace operator.

Factional Laplace operators appear specially in Finance, Physics and Ecology, see for example \cite{A}. These operators have attracted special attention during the last decade and a rich theory has been developed after a breakthrough was given when Caffarelli and Silvestre \cite{CaSi} introduced a characterization of the fractional Laplace operator $(-\Delta)^{s}$ in terms of a Dirichlet-to-Neumann map associated to a suitable extension problem for $0 < s < 1$.

Simultaneously, a lot of attention also has been paid to nonlinear problems of the form
\begin{equation}\label{3}
\left\{
\begin{array}{rrll}
(-\Delta)^{s} u &=& f(x,u) & {\rm in} \ \ \Omega\\
u &=& 0  & {\rm in} \ \ \R^n\setminus\Omega
\end{array}
\right.
\end{equation}
For works on existence of solutions we refer to \cite{ros131, ros137, ros138, auto2, auto1, ros215, niang, ros264, ros265, ros267, val}, nonexistence to \cite{ros129, val},  symmetry to \cite{ros21, ros95}, regularity to \cite{ros10, cabre1, ROS}, and other qualitative properties of solutions to \cite{ros130, ros1}. For developments related to (\ref{3}) involving the spectral fractional Laplace operator $\mathcal{A}^{s}$, we refer to \cite{BaCPS, BrCPS, CT, CDDS, choi, CK, tan1, T, yang}, among others.

Systems like (\ref{1}) are strongly coupled vector counterparts closely related to (\ref{3}) which have been studied by several authors during the two last decades for $s = t = 1$ (we refer to the survey \cite{DG} and references therein). More specifically, a priori bounds and existence of positive solutions have been addressed in these cases. In view of what is known for scalar equations and for systems of the type (\ref{1}) with $s = t = 1$, one expects that a priori bounds depend on the values of the exponents $p$ and $q$. Indeed, the values $p$ and $q$ should be related to Sobolev embedding theorems.

A rather classical fact is that a priori bounds allow to establish existence of positive solutions for systems by mean of topological methods such as degree theory and Krasnoselskii's index theory. For a list of works concerning with non-variational elliptic systems involving Laplace operators we refer to \cite{BM, CFM, DGS1, DGS2, M, sirakov, Sira, zou} and references therein.

One of the goals of this work is to establish existence of positive classical solutions of non-variational elliptic systems of the type (\ref{1}) by mean of a priori bounds for a family of exponents $p$ and $q$. By a classical solution of the system (\ref{1}), we mean a couple $(u,v) \in (C^\alpha(\R^n))^2$, with $0<\alpha<1$, satisfying (\ref{1}) in the usual sense.

Our main result is

\begin{teor}\label{teo1}
Let $\Omega$ be a bounded open subset of $C^2$ class of $\R^n$. Assume that $n \geq 2$, $0 < s, t < 1$, $n > 2s + 1$, $n > 2t + 1$, $p,q \geq 1$, $pq > 1$ and either
\begin{equation} \label{4}
\left(\frac{2s}{p}+2t\right)\frac{p}{pq-1}\geq n-2s\ \ {\rm or}\ \ \left(\frac{2t}{q}+2s\right)\frac{q}{pq-1}\geq n-2t\, .
\end{equation}
Then, the system (\ref{1}) admits, at least, one positive classical solution. Moreover, all such solutions are uniformly bounded in the $L^{\infty}$-norm by a constant that depends only on $s, t, p, q$ and $\Omega$.
\end{teor}

\begin{rem}
Systems like (\ref{1}) have been widely studied when $s = t = 1$. In this case, it arises the well-known notions of superlinearity and criticality (under the form of critical hyperbole), see \cite{Mi1, Mi2, SZ1}. In particular, one knows that (see \cite{DF, DGR, vander, HMV}) the system (\ref{1}) always admits a positive classical solution provided that $pq > 1$ and

\[
\frac{1}{p+1} + \frac{1}{q+1} > \frac{n - 2}{n}\, .
\]
\end{rem}

\begin{rem}
When $0 < s = t < 1$ and $p, q > 1$, a priori bounds and existence of positive classical solutions of (\ref{1}) for the spectral fractional Laplace operator $\mathcal{A}^s$ have been derived in \cite{choi} provided that

\begin{equation} \label{Hip}
\frac{1}{p+1} + \frac{1}{q+1} > \frac{n - 2s}{n}\, .
\end{equation}
\end{rem}

\begin{rem}
When $0 < s = t < 1$ and $pq > 1$, the condition (\ref{4}) implies (\ref{Hip}).
\end{rem}

The approach used in the proof of Theorem \ref{teo1} is based on the blow-up method, firstly introduced by Gidas and Spruck in \cite{GS1} to treat the scalar case and later extended to strongly coupled systems like (\ref{1}) with $s=t=1$ in \cite{M} and then in \cite{DGS1, DGS2, Sira, zou}.  The method is used to get uniform bounds of solutions through a contradiction argument by assuming that desired bounds fail and relies on Liouville theorems for related problems in the whole space $\R^n$ and in half-spaces. The proof of these results is usually the most difficult part in applying the blow-up method.

For this purpose, we first shall establish Liouville theorems for the system

\begin{equation}  \label{sistema2}
\left\{\begin{array}{llll}
(-\Delta)^{s} u = v^{p} & {\rm in} \ \ G\\
(-\Delta)^{t} v = u^{q} & {\rm in} \ \ G\\
\end{array}\right.
\end{equation}
for $G=\R^{n}$ and $G=\R^{n}_{+} = \{ x =(x_1, \ldots, x_n) \in \R^n:\ x_n > 0\}$. In that latter, we assume the Dirichlet condition $u = 0 = v$ in $\R^n\setminus\R^{n}_{+}$.

We recall that a viscosity super-solution of the above system is a couple $(u,v)$ of continuous functions in $\R^n$ such that $u,v \geq 0$ in $\R^{n} \setminus \overline{G}$ and for each point $x_{0}\in G$ there exists a neighborhood $U$ of $x_{0}$ with $\overline{U} \subset G$ such that for any $\varphi,\psi \in C^{2}(\overline{U})$ satisfying $u(x_{0}) = \varphi(x_{0})$, $v(x_{0}) = \psi(x_{0})$, $u \geq \varphi$ and $v \geq \psi$ in $U$, the functions defined by

\begin{equation} \label{2.2}
\overline{u}=\left\{\begin{array}{lllc}
\varphi & {\rm in} \ U  \\
u & {\rm in} \ \R^{n}\setminus U
\end{array}\right. \ \ {\rm and}\ \
\overline{v}=\left\{\begin{array}{lllc}
\psi & {\rm in} \ U  \\
v & {\rm in} \ \R^{n}\setminus U
\end{array}\right.
\end{equation}

\n satisfy
\[
(-\Delta)^{s}\overline{u}(x_{0}) \geq v^{p}(x_{0}) \ \ {\rm and}\ \
(-\Delta)^{t}\overline{v}(x_{0}) \geq u^{q}(x_{0})\, .
\]

\n In a natural way, we have the notions of viscosity sub-solution and viscosity solution.

In the following, we state two Liouville theorems for the system (\ref{sistema2}).

\begin{teor}\label{Teo2}
Assume that $n \geq 2$, $0 < s,t < 1$, $n > 2s$, $n > 2t$, $p,q > 0$ and $pq > 1$. Then, the only non-negative viscosity super-solution of the system (\ref{sistema2}) with $G=\R^{n}$ is the trivial if and only if (\ref{4}) holds.
\end{teor}

\begin{teor}\label{Teo3}
Assume that $n \geq 2$, $0 < s,t < 1$, $n > 2s + 1$, $n > 2t + 1$, $p,q \geq 1$ and $pq > 1$. If the condition (\ref{4}) holds, then the only non-negative viscosity bounded solution of the system (\ref{sistema2}) with $G=\R^{n}_{+}$ is the trivial.
\end{teor}

\begin{rem}
Non-existence results of positive solutions have been established for the scalar problem

\[
(-\Delta)^{s} u = u^p \ \ {\rm in} \ \ G
\]

\n in both cases $G = \R^n$ and $G = \R^n_{+}$ by assuming that $n > 2s$ and $1 < p < \frac{n + 2s}{n - 2s}$, see \cite{GS, GS1} for $s=1$ and \cite{JLX, xia} for $0 < s < 1$.
\end{rem}

\begin{rem}
A number of works has focused on non-existence of positive solutions of (\ref{sistema2}) for $G = \R^n$ and $G = \R^n_{+}$ when $s = t = 1$ and $0 < s = t < 1$. We refer for instance to \cite{BM, BuMa, DF1, Mi2, PQS, SZ2, Sou, So} for $s = t = 1$ and \cite{xia} for $0 < s = t < 1$, among other references.
\end{rem}

Several methods have been employed in the proof of non-existence results of positive solutions of elliptic systems. Our approach is inspired on a method developed by Quaas and Sirakov in \cite{sirakov} to treat systems involving different uniformly elliptic linear operators based on maximum principles. Particularly, some maximum principles for fractional operators proved by Silvestre \cite{S} and Quaas and Xia \cite{xia} as well as some fundamental lemmas to be proved in the next section will play an important role in the proof of Theorem \ref{Teo2} and also in the remainder of the work.

The paper is organized into four sections. In Section 2 we prove three lemmas required in the proof of Theorem \ref{Teo2} which will be presented in Section 3. Section 4 is devoted to the proof of Theorem \ref{Teo3}. In Section 5, we prove Theorem \ref{teo1} by using Theorems \ref{Teo2} and \ref{Teo3}.\\

\section{Fundamental lemmas}

We next present three lemmas which will be used in the proof of Theorem \ref{Teo2}.

Throughout the paper, it is assumed that $p, q > 0$ and $pq > 1$. So, thanks to a suitable rescaling of $u$ and $v$, we can assume that $C(n,s) = 1$ and $C(n,t) = 1$.

Given a non-negative continuous function $u : \R^n \rightarrow \R$, define

\[
m_u(r) = \displaystyle \min_{\vert x\vert\leq r}u(x)\,
\]
for $r>0$.
\begin{lema}\label{lema4.1}
Let $0 < s < 1$, $n>2s$ and $u\neq 0$ be a non-negative viscosity super-solution of
\begin{equation}\label{equacao}
(-\Delta)^{s} u = 0\ {\rm in}\ \R^{n}\, .
\end{equation}
Then, for each $R_0 > 1$ and $\sigma \in (-n,-n+2s)$, there exists a constant $C > 0$, independent of $u$, such that
\begin{equation}
m_u(r) \geq C m_u(R_{0}) r^{\sigma} \label{16}
\end{equation}

\n for all $r \geq R_0$.
\end{lema}

\n By a non-negative viscosity super-solution of the equation (\ref{equacao}), we mean a non-negative continuous function $u: \R^n \rightarrow \R$ satisfying the following property: each point $x_{0} \in \R^n$ admits a neighborhood $U$ such that for any function $\varphi \in C^{2}(\overline{U})$ with $u(x_{0}) = \varphi(x_{0})$ and $u \geq \varphi$ in $U$, the function defined by

\[
\overline{u}=\left\{\begin{array}{lllc}
\varphi & {\rm in} \ U  \\
u & {\rm in}  \ \R^{n}\setminus U
\end{array}\right.
\]

\n satisfies
\[
(-\Delta)^{s}\overline{u}(x_{0}) \geq 0\, .
\]

\n {\bf Proof of Lemma \ref{lema4.1}.} Let $R_0$, $\sigma$ and $u$ be as in the above statement. Given $R > R_0$ and $\varepsilon > 0$, we consider the function
\begin{equation}
w(r)=\left\{\begin{array}{lc}
\varepsilon^{\sigma} $   if   $ 0 < r \leq \varepsilon & \\
r^{\sigma}$   if   $ \varepsilon \leq r
\end{array}\right.
\end{equation}
We first assert that $(-\Delta)^{s} w(r) < 0$ for all $R_0 < r < R$ and $\varepsilon > 0$ small enough. In fact, for $|x|=r$, we have

\begin{eqnarray*}
2(-\Delta)^{s} w(r)&=& - \int\limits_{B_{\varepsilon}(-x)}\frac{\varepsilon^{\sigma}}{|y|^{n+2s}}\; dy - \int\limits_{B_{\varepsilon}(x)}\frac{\varepsilon^{\sigma}}{|y|^{n+2s}}\; dy - \int\limits_{B_{\varepsilon}^{c}(-x)}\frac{|x+y|^{\sigma}}{|y|^{n+2s}}\; dy\\
&& - \int\limits_{B_{\varepsilon}^{c}(x)}\frac{|x-y|^{\sigma}}{|y|^{n+2s}}\; dy + 2\int\limits_{\R^{n}}\frac{|x|^{\sigma}}{|y|^{n+2s}}\; dy\\
&=& - \int\limits_{\R^{n}}\frac{|x+y|^{\sigma} + |x-y|^{\sigma} -2 |x|^{\sigma}}{|y|^{n+2s}}\; dy \\
&& + \left(\int\limits_{B_{\varepsilon}(-x)}\frac{|x+y|^{\sigma} -\varepsilon^{\sigma}}{|y|^{n+2s}}\; dy + \int\limits_{B_{\varepsilon}(x)}\frac{|x-y|^{\sigma} -\varepsilon^{\sigma}}{|y|^{n+2s}}\; dy\right) \\
&=& 2(-\Delta)^{s} |x|^{\sigma} + \left(\int\limits_{B_{\varepsilon}(-x)}\frac{|x+y|^{\sigma} -\varepsilon^{\sigma}}{|y|^{n+2s}}\; dy + \int\limits_{B_{\varepsilon}(x)}\frac{|x-y|^{\sigma} -\varepsilon^{\sigma}}{|y|^{n+2s}}\; dy\right)\, .
\end{eqnarray*}

\n Since  $R_0 > 1$, the two last above integral converge uniformly to $0$ for $|x| > R_{0}$ as $\varepsilon \rightarrow 0$.

On the other hand, using that $R_0 > 1$ and $\sigma \in (-n,-n+2s)$ and the fact that $|x|^{-n+2s}$ is the fundamental solution of the fractional Laplace operator $(-\Delta)^{s}$ (see \cite{CaSi}), one easily checks that $(-\Delta)^{s} |x|^{\sigma} < 0$ for all $|x| > R_0$, see \cite{FQ}. Thus, the above claim follows for $\varepsilon > 0$ small enough.

For such a parameter $\varepsilon$ and $|x|=r$, we set

\[
\varphi(x) = m_u(R_{0}) \frac{w(r) - w(R)}{w(\varepsilon) - w(R)}
\]

\n for all $|x| < R$ and $\varphi(x)=0$ for $|x| \geq R$. As can easily be checked, $(-\Delta)^{s} \varphi \leq 0$ for all $R_0 < |x| < R$. Moreover, we have $u(x) \geq \varphi(x)$ for $|x| \leq R_0$ or $|x| \geq R$, so that the Silvestre's strong maximum principle \cite{S} readily yields $u(x) \geq \varphi(x)$ for all $R_0 \leq |x| \leq R$. Finally, letting $R \rightarrow \infty$ in this last inequality, we achieve the expected conclusion with $C = \varepsilon^{-\sigma}$. \fim\\

Our second auxiliary lemma is

\begin{lema}\label{lema4.2}
Let $0 < s < 1$, $n>2s$ and $u\neq 0$ be a non-negative viscosity super-solution of (\ref{equacao}). Then, there exist constants $C>0$ and $R_0 > 0$, independent of $u$, such that
\begin{equation} \label{17}
m_u(r/2) \leq C m_u(r)
\end{equation}
\n for all $r \geq R_0$.
\end{lema}
\n {\bf Proof of Lemma \ref{lema4.2}.} Given $r > 0$ and $\varepsilon>0$, set

\[
R=r\left[\frac{\varepsilon}{1+\varepsilon 2^{-n+2s}}\right]^{1/(n-2s)}\, ,
\]

\n where $\varepsilon$ is chosen such that $R < r/2$.

Consider the functions

\[
w_{r}(\overline{r}) = \left\{
\begin{array}[c]{ll}%
(R)^{-n+2s} & {\rm if}\ 0<\overline{r}\leq R\\
\overline{r}^{-n+2s} & {\rm if}\ R\leq \overline{r}\leq 2r\\
(2r)^{-n+2s} & {\rm if} \ \overline{r}\geq 2r
\end{array}
\right.
\]
and

\[
w(\overline{r}) = \left\{
\begin{array}[c]{ll}%
(R)^{-n+2s} & {\rm if}\ 0<\overline{r}\leq R\\
\overline{r}^{-n+2s} & {\rm if}\ R\leq \overline{r}
\end{array}
\right.
\]

\n Given a fixed function $u$ as in the above statement, we define

\[
\varphi(x) = m_u(r/2) \frac{w_{r}(\overline{r}) - w(2r)}{w(R) - w(2r)}
\]

\n for $x$ with $|x| = \overline{r}$. As a direct consequence, one has $u(x) \geq \varphi(x)$ for all $x$ with $\vert x\vert \leq r/2$ and $\vert x\vert \geq 2r$. Moreover, decreasing $\varepsilon$, if necessary, one gets

\begin{eqnarray*}
2(-\Delta)^{s} w_{r}(\overline{r})= - \int\limits_{B_{R}(x)} \frac{r^{-n+2s}}{\varepsilon |y|^{n+2s}}+\frac{(2r)^{-n+2s}}{|y|^{n+2s}}\; dy - \int\limits_{B_{2r}^{c}(x)}\frac{(2r)^{-n+2s}}{|y|^{n+2s}}\; dy \\
- \int\limits_{B_{2r}(x)\setminus B_{R}(x)}\frac{| x-y|^{-n+2s}}{|y|^{n+2s}}\; dy + \int\limits_{\R^{n}}\frac{\vert x\vert^{-n+2s}}{|y|^{n+2s}}\; dy \leq 0
\end{eqnarray*}
for all $r/2 < \overline{r} < 2r$. Thus, $(-\Delta)^{s}\varphi(x) \leq 0$ for all $x$ with $r/2 < \vert x\vert < 2r$.

Evoking the Silvestre's maximum principle \cite{S}, we then deduce that $u(x) \geq \varphi(x)$ for all $x$ with $r/2 < \vert x\vert < 2r$. Lastly, we assert that this conclusion leads to

\[
m_u(r) \geq \varepsilon m_u(r/2)(1-2^{-n+2s})\, .
\]

\n In fact, we have

\[
\varphi(x) = m_u(r/2) \geq \varepsilon m_u(r/2)(1-2^{-n+2s})
\]

\n if $0 < \vert x\vert \leq R$, and

\[
\varphi(x) = \varepsilon m_u(r/2) \frac{\overline{r}^{-n+2s}-(2r)^{-n+2s}}{r^{-n+2s}} \geq \varepsilon m_u(r/2)(1-2^{-n+2s})
\]

\n if $R < \vert x\vert \leq r$. So, the result follows with $C = (\varepsilon(1-2^{-n+2s}))^{-1}$ by minimizing $u$ on the closed ball $|x| \leq r$.\; \fim\\

Our third lemma concerns with the behavior of fractional Laplace operators applied to the function $\Theta(x) = \log(1+ |x|) |x|^{-n+2s}$.

\begin{lema}\label{lema6.1.}
Let $0 < s < 1$ and $n>2s$. Then, there exists a constant $C_0 > 0$ such that
\[
(-\Delta)^{s}\Theta (x) \leq C_0 |x|^{-n}
\]

\n for all $x \neq 0$.
\end{lema}

\n {\bf Proof of Lemma \ref{lema6.1.}.} Using that $|x|^{-n + 2s}$ is the fundamental solution of $(-\Delta)^{s}$ (see \cite{CaSi}), one first has
\begin{eqnarray*}
-2(-\Delta)^{s} \Theta(x) &=& \int\limits_{\R^{n}}\frac{\log(1+ |x-y|) |x-y|^{-n+2s}}{|y|^{n+2s}}\; dy\\
&& + \int\limits_{\R^{n}}\frac{\log(1+ |x+y|) |x+y|^{-n+2s}}{|y|^{n+2s}}\; dy - 2\int\limits_{\R^{n}}\frac{\log(1+|x|) |x|^{-n+2s}}{|y|^{n+2s}}\; dy\\
&= & \int\limits_{\R^{n}}\frac{\left(\log(1+ |x-y|)-\log(1+ |x|)\right) |x-y|^{-n+2s}}{|y|^{n+2s}}\; dy\\
&& + \int\limits_{\R^{n}}\frac{\left(\log(1+ |x+y|)-\log(1+ |x|)\right) |x+y|^{-n+2s}}{|y|^{n+2s}}\; dy\\
&=& \int\limits_{\R^{n}}\left(\log\left(\frac{1+ |x-y|}{1+ |x|}\right) |x-y|^{-n+2s}\right)\frac{1}{|y|^{n+2s}}\; dy\\
&& + \int\limits_{\R^{n}}\left(\log\left(\frac{1+ |x+y|}{1+ |x|}\right) |x+y|^{-n+2s}\right)\frac{1}{|y|^{n+2s}}\; dy\\
&=& \int\limits_{\R^{n}}r^{-n}\left(\log\left(\frac{1+r |e_{1}-z|}{1+ r}\right) |e_{1}-z|^{-n+2s}\right)\frac{1}{|z|^{n+2s}}\; dz\\
&& + \int\limits_{\R^{n}}r^{-n}\left(\log\left(\frac{1+r |e_{1}+z|}{1+ r}\right) |e_{1}+z|^{-n+2s}\right)\frac{1}{|z|^{n+2s}}\; dz\, ,
\end{eqnarray*}
where $x=re_{1}$ and $z=y/r.$ Note that there is no loss of generality in considering $x=re_{1}$, since $\log(1+ |x|)$ and $|x|^{-n+2s}$ are radially symmetric.

In order to complete the proof we just need to find a constant $C_0 > 0$ such that

\begin{equation} \label{6.8}
\int\limits_{\R^{n}}\frac{\left(\log\left(\frac{1+r |e_{1}-z|}{1+ r}\right) |e_{1}-z|^{-n+2s} + \log\left(\frac{1+r |e_{1}+z|}{1+ r}\right) |e_{1}+z|^{-n+2s}\right)}{|z|^{n+2s}}\; dz \geq -C_0\, .
\end{equation}

\n For this purpose, we write for $\rho > 0$, $\gamma \in [0,1)$ and $r \geq 0$,

\begin{equation}
\log\left(\frac{1+r |e_{1}-z|}{1+ r}\right) |e_{1}-z|^{-n+2s} = g(|e_{1}-z|,\gamma)
\end{equation}

\n and

\begin{equation}
\log\left(\frac{1+r |e_{1}+z|}{1+ r}\right) |e_{1}+z|^{-n+2s}=g(|e_{1}+z|,\gamma)\, ,
\end{equation}

\n where

\[
g(\rho,\gamma)=\rho^{-n+2s}\log(1+\gamma(\rho-1))
\]

\n and

\[
\gamma = \frac{r}{1+r}\, .
\]

\n Consider first $B_{1}=\{z : |z+e_{1}| \leq 1/2\}$ and note that $g(|e_{1}-z|,\gamma)$ is bounded in $B_{1}$, while $g(|e_{1}+z|,\gamma)$ has a singularity at $-e_{1} \in B_{1}$. Then, for some constants $C>0$, independent of $\gamma$, we have

\begin{eqnarray*}
\int\limits_{B_{1}}\frac{| g(|e_{1}+z|,\gamma)|}{\vert z\vert^{n+2s}}\; dz&=&
\int\limits_{B_{1/2}(0)}\frac{| g(|z| ,\gamma)|}{|z-e_{1}|^{n+2s}}\; dz \leq - C \int\limits_{0}^{1/2}g(\rho,\gamma)\rho^{n-1}\; d\rho\\
&\leq & -C \int\limits_{0}^{1/2}\rho^{2s-1}\log(\rho)\; d\rho \leq C\, .
\end{eqnarray*}

\n Since $1 + \gamma (\rho-1) \geq \rho$ as $\gamma \in [0,1)$, the integral in (\ref{6.8}), when considered over $B_{1}$, is bounded below by a constant independent of $r$. In a similar way, the conclusion follows for the set $B_{2} = \{z: |z-e_{1}| \leq 1/2\}$.

On the set $B_{3} =\{z: |z| \geq 2\}$, for some constant $C>0$, independent of $\gamma$, we have

\[
|g(|e_{1}-z|,\gamma) +g(|e_{1}+z|,\gamma)| \leq C |z|^{-2n} \log(|z|)\, .
\]
Thus, the integral in (\ref{6.8}), when considered over $B_{3}$, is also bounded below by a constant independent of $r$.

It then remains to analyze the behavior of the integral over $B_{4}=\{z : |z| \leq 1/2 \}$.

\n For each fixed $r \geq 0$ and $\gamma \in [0,1)$,  define $f_{r}:\R^{n} \rightarrow \R$ given by $f_{r}(z)= g(|e_{1}+z|, \gamma) + g(|e_{1}-z|, \gamma)$. Using that $f_{r}(0)=0$ and $D(f_{r}(0))=0$, the Taylor formula provides

\begin{equation}
f_{r}(z)=z^{t} \cdot \int\limits_{0}^{1}(1-\rho)D^{2}(f_{r}(\rho z))\; d\rho\cdot z\, ,
\end{equation}

\n where all derivatives are taken only with respect to the variable $z$. Thus, the estimate of the integral (\ref{6.8}) over $B_{4}$ follows if we can show that

\begin{equation} \label{6.12}
\left|\frac{\partial^{2}f_{r}(z)}{\partial z_{i}\partial z_{j}}\right| \leq C
\end{equation}

\n for all $|z| \leq 1/2$, where $C >0$ is a constant independent of $r$.

On the other hand, a straightforward computation gives

\[
\frac{d}{d\rho}g(\rho,\gamma) = (-n+2s)\rho^{-n+2s-1}\log(1+\gamma(\rho-1))+\frac{\gamma\rho^{-n+2s}}{1+\gamma(\rho-1)}
\]

\n and

\begin{eqnarray*}
\frac{d^2}{d\rho^2}g(\rho,\gamma) &=& (-n+2s)(-n+2s-1)\rho^{-n+2s-2}\log(1+\gamma(\rho-1))\\
&& + \frac{2\gamma (-n+2s)\rho^{-n+2s-1}}{1+\gamma(\rho-1)}-\frac{\gamma^{2}\rho^{-n+2s}}{(1+\gamma(\rho-1))^{2}}\, .
\end{eqnarray*}

\n Then, one easily checks that

\[
| \frac{d}{d\rho}g(\rho,\gamma)|,\ |\frac{d^2}{d\rho^2}g(\rho,\gamma)|\leq C
\]

\n for all $1/2 \leq \rho \leq 3/2$ and $\gamma \in [0,1)$, where $C$ is a constant independent of $\rho$ and $\gamma$.

\n So, for certain bounded functions $D_{ij}$ and $d_{ij}$ in $B_{4}$, we have

\[
\frac{\partial^{2}f_{r}(z)}{\partial z_{i}\partial z_{j}} = \frac{d^2}{d\rho^2}g( |e_{1}+z|,\gamma)D_{ij} + \frac{d}{d\rho}g( |e_{1}+z|,\gamma)d_{ij}
\]

\n and (\ref{6.12}) follows.

Finally, joining the above estimates on the four sets $B_i$, one gets (\ref{6.8}) as desired.\; \fim\\

\section{Proof of Theorem \ref{Teo2}}

We organize the proof of Theorem \ref{Teo2} into two stages, according to the sufficiency and necessity of the assumption (\ref{4}).

\n {\bf Proof of the sufficiency of (\ref{4}).} We analyze separately two different cases:

\begin{itemize}
\item[(I)] $\left(\frac{2s}{p}+2t\right)\frac{p}{pq-1}>n-2s$ or $\left(\frac{2t}{q}+2s\right)\frac{q}{pq-1}>n-2t$;

\item[(II)] $\left(\frac{2s}{p}+2t\right)\frac{p}{pq-1} = n-2s$ or $\left(\frac{2t}{q}+2s\right)\frac{q}{pq-1} = n-2t$.
\end{itemize}

We first assume the situation (I). Let $(u,v)$ be a non-negative viscosity super-solution of the system (\ref{sistema2}) with $G=\R^{n}$ and $\eta :[0,+\infty)\rightarrow \R$ be a $C^\infty$ cutoff function satisfying $0 \leq \eta \leq 1$, $\eta$ is non-increasing, $\eta(r) = 1$ if $0 \leq r \leq 1/2$ and $\eta(r)=0$ if $r\geq 1$. Clearly, there exists a constant $C > 0$ such that $(-\Delta)^{s} \eta( |x| )\leq C$ and $(-\Delta)^{t} \eta( |x| )\leq C$.

Choose $R_0 > 0$ as in Lemma \ref{lema4.2} for $s$ and $t$, simultaneously, and consider the functions

\[
\xi_{u}(x) = m_{u}(R_0/2)\eta(|x|/R_0) \ {\rm and}\ \xi_{v}(x) = m_{v}(R_0/2)\eta(|x|/R_0)\, .
\]

\n For some constant $C_0 > 0$, independent of $R_0$, $u$ and $v$, we have

\[
(-\Delta)^{s}(\xi_{u}(x)) \leq C_0 \frac{m_{u}(R_0/2)}{{R_0}^{2s}}\ {\rm and}\  (-\Delta)^{t}(\xi_{v}(x)) \leq C_0 \frac{m_{v}(R_0/2)}{{R_0}^{2t}}\, .
\]
Moreover, $\xi_{u}(x)=0\leq u(x)$ if $|x| > R_0$ and $\xi_{u}(x) = m_{u}(R_0/2) \leq u(x)$ if $|x| \leq R_0/2$. Similarly, $\xi_{v}(x)=0 \leq v(x)$ if $|x| > R_0$ and $\xi_{v}(x) = m_{v}(R_0/2) \leq v(x)$ if $|x| \leq R_0/2$. Thus, the functions $u-\xi_{u}$ and $v-\xi_{v}$ attain their global minimum values at points $x_{u}$ and $x_{v}$ with $|x_{u}| < R_0$ and $|x_{v}| < R_0$, respectively.

Now let $\varphi(x):= \xi_{u}(x) - \xi_{u}(x_{u}) + u(x_{u})$ and $\psi(x):= \xi_{v}(x) -\xi_{v}(x_{v}) + v(x_{v})$. Note that $\varphi(x_{u}) = u(x_{u})$, $\psi(x_{v}) = v(x_{v}), u(x) \geq \varphi(x)$ and $v(x) \geq \psi(x)$ for all $x \in B(0,R_0)$. Let $\overline{u}$ and $\overline{v}$ be defined as in (\ref{2.2}) with $U = B(0,R_0)$. Since $(u,v)$ is a viscosity super-solution of (\ref{sistema2}), one has

\begin{equation} \label{6.1}
(-\Delta)^{s} (\overline{u})(x_{u}) \geq v^{p}(x_{u})\ {\rm and}\ (-\Delta)^{t} (\overline{v})(x_{v}) \geq u^{q}(x_{v})\, .
\end{equation}

We now assert that

\[
(-\Delta)^{s} (\overline{u})(x_{u}) \leq (-\Delta)^{s}(\xi_{u})(x_{u})\ {\rm and}\ (-\Delta)^{t} (\overline{v})(x_{v}) \leq (-\Delta)^{t}(\xi_{v})(x_{v})\, .
\]
In fact, note that $w_{u}(x) := \overline{u}(x)-\xi_{u}(x) \geq 0$ for all $x\in \R^{n}$ and $x_{u}$ is a global minimum point of $w_{u}$. Thus, we have $(-\Delta)^{s}(w_{u})(x_{u}) \leq 0$ and thus the first inequality follows. The other inequality also follows in an analogous way. Therefore, from (\ref{6.1}), one gets

\begin{equation}\label{16}
m_{u}^{q}(R_0) \leq u^{q}(x_{v}) \leq C_0 \frac{m_{v}(R_0/2)}{{R_0}^{2t}}\ {\rm and}\ m_{v}^{p}(R_0) \leq v^{p}(x_{u}) \leq C_0 \frac{m_{u}(R_0/2)}{{R_0}^{2s}}\, .
\end{equation}

\n Applying Lemma \ref{lema4.2} in the above inequalities, one then derives
\begin{equation} \label{6.2}
m_{u}(R_0) \leq \frac{C_{1}}{{R_0}^{\left(\frac{2s}{p}+2t\right)\frac{p}{pq-1}}} \ {\rm and}\ m_{v}(R_0) \leq \frac{C_{2}}{{R_0}^{\left(\frac{2t}{q}+2s\right)\frac{q}{pq-1}}}\, .
\end{equation}

We now consider the case (I). It suffices to assume that $(\frac{2s}{p} + 2t) \frac{p}{pq-1} > n-2s$, since the argument is analogous for the second inequality in (I). Choose $-n < \sigma_{1} < -n+2s$ such that
\[
\left(\frac{2s}{p}+2t\right)\frac{p}{pq-1}+\sigma_{1} > 0\, .
\]

\n By Lemma \ref{lema4.1}, we have
\[
m_{u}(r)\leq m_{u}(R_0)\leq\frac{C}{r^{\left(\frac{2s}{p}+2t\right)\frac{p}{pq-1}+\sigma_{1}}}
\]
for all $r \geq R_0 \geq 1$. Therefore, $m_{u}(r)$ goes to $0$ as $r \rightarrow + \infty$, providing the contradiction $(u,v)=(0,0)$.

Finally, assume the situation (II). In a similar way, we analyze only the equality $(\frac{2s}{p}+2t) \frac{p}{pq-1} = n-2s$. Let $(u,v)$ be non-negative viscosity super-solution of (\ref{sistema2}) with $G=\R^{n}$. We begin by proving that for certain $C>0$ and $R_0 > 0$, we have

\begin{equation}\label{15}
m_{u}(r)\geq C m_{u}(R_0)r^{-n+2s}
\end{equation}
for all $r \geq R_0$. Indeed, by Lemma \ref{lema4.1} and (\ref{16}), for any $-n < \sigma < -n+2s$, we have

\begin{equation}\label{14.}
(-\Delta)^{s}u(x) \geq v^{p}(x) \geq m_{v}(r)^{p} \geq C(m_{u}(2r))^{pq}r^{2tp} \geq C(m_{u}(R_0))^{pq}r^{\sigma pq+2tp}
\end{equation}
for all $x$ with $|x| = r \geq R_0$.

Now consider the function

\begin{equation}
w(r)=\left\{\begin{array}{lc}
\varepsilon^{-n+2s} $   if   $ 0 < r \leq \varepsilon & \\
r^{-n+2s}$   if   $ \varepsilon \leq r
\end{array}\right.
\end{equation}
where $0 < \varepsilon < R_0/2$. Since $|x|^{-n+2s}$ is the fundamental solution of the fractional Laplace operator $(-\Delta)^{s}$ (see \cite{CaSi}), we have

\begin{eqnarray*}
2(-\Delta)^{s} w(r)= \left(\int\limits_{B_{\varepsilon}(-x)}\frac{|x+y|^{-n+2s} -\varepsilon^{-n+2s}}{|y|^{n+2s}}\; dy + \int\limits_{B_{\varepsilon}(x)}\frac{|x-y|^{-n+2s} -\varepsilon^{-n+2s}}{|y|^{n+2s}}\; dy\right)\, ,
\end{eqnarray*}
where $|x| = r$. It is clear that $|y| \geq |x|/2$ whenever $|x| \geq R_0$ and $y \in B_{\varepsilon}(x)$. Thus,
\[
\int\limits_{B_{\varepsilon}(x)}\frac{|x-y|^{-n+2s} -\varepsilon^{-n+2s}}{|y|^{n+2s}}\; dy \leq \frac{C}{r^{n+2s}}
\]
for some constant $C > 0$ and then, by symmetry of the integrals, one obtains

\[
2(-\Delta)^{s} w(r)\leq\frac{C}{r^{n+2s}}\, .
\]
For fixed $R_1 > R_0$, we define the comparison function

\[
\varphi(x) = m_u(R_0) \frac{w(r) - w(R_1)}{w(\varepsilon) - w(R_1)}
\]

\n for all $x$ with $|x| < R_1$ and $\varphi(x)=0$ for $|x| \geq R_1$. As can easily be checked,
\begin{equation}\label{13}
(-\Delta)^{s} \varphi(x) \leq\frac{C_{1}}{|x|^{n+2s}}
\end{equation}
for all $x$ with $R_0 < |x| < R_1$. On the other hand, since $n = pq(n-2s) - 2tp$, we can choose $\sigma \in (-n,-n+2s)$ such that $-\sigma pq-2tp<n+2s$. Then, using (\ref{14.}) and (\ref{13}), one gets
\[
(-\Delta)^{s} \varphi(x) \leq\frac{C_{1}}{\vert x\vert^{n+2s}}\leq\frac{C_{1}}{\vert x\vert^{-\sigma pq-2tp}}\leq (-\Delta)^{s}u(x)
\]
for all $x$ with $R_0 < |x| < R_1$ and $u(x) \geq \varphi(x)$ for $|x| \leq R_0$ or $|x| \geq R_1$, so that the Silvestre's maximum principle \cite{S} readily yields $u(x) \geq \varphi(x)$ for all $R_0 \leq |x| \leq R_1$. Finally, letting $R_1 \rightarrow + \infty$ in this last inequality, the claim (\ref{15}) follows.

In the sequel, we split the proof into two cases according to the value of $-n+2s$. The first one corresponds to $-n+2s \in (-n,-1]$. In this range, note that the function $\Theta$, defined above Lemma \ref{lema6.1.}, is decreasing for all $r>0$, with a singularity at the origin if $-n+2s\in(-n,-1)$ and bounded if $-n+2s=-1$. For $0 < \varepsilon < R_0/2$, we define the function
\[
w(r) = \left\{
\begin{array}[c]{ll}%
\Theta(\varepsilon) & {\rm if}\ 0<r\leq\varepsilon\\
\Theta(r) & {\rm if}\ \varepsilon < r
\end{array}
\right..
\]
Using Lemma \ref{lema6.1.}, for any $r \geq R_0$ and $x$ with $|x| = r$, we have
\begin{eqnarray*}
(-\Delta)^{s} w(r) & \leq & \int\limits_{B_{\varepsilon}(x)}\frac{\log(1+ |x-y|) |x-y|^{-n+2s} - \log(1+\varepsilon)\varepsilon^{-n+2s}}{ |y|^{n+2s}}\; dy + \frac{C}{r^n}\\
& \leq & C \frac{\varepsilon^{2s}}{r^{n+2s}} + \frac{C}{r^n} \leq \frac{C}{r^n}
\end{eqnarray*}
\n for all $r \geq R_0$ and some constant $C>0$ independent of $r$.

Let $\varphi$ be defined as above for $R_1 > R_0$. Again, we have $\varphi(x) \leq u(x)$ for all $x$ with $|x| \leq R_0$ or $|x| \geq R_1$. Moreover,

\begin{equation} \label{24}
(-\Delta)^{s} \varphi(x) \leq \frac{C}{|x|^{n}}
\end{equation}
for all $x$ with $R_0 < |x| < R_1$. From (\ref{14.}), one also has

\begin{equation}\label{14}
(-\Delta)^{s}u(x) \geq C(m_{u}(R_0))^{pq}r^{(-n+2s)pq+2tp}=\frac{C}{|x|^{n}}
\end{equation}
for $r \geq R_0$. By Silvestre's maximum principle \cite{S}, we derive $u(x) \geq \varphi(x)$ for all $R_0 < |x| < R_1$. Letting $R_1 \rightarrow +\infty$ in this inequality, one obtains

\[
u(x) \geq C \frac{\log(1+ |x|)}{|x|^{n-2s}}\, .
\]

\n On the other hand, using (\ref{6.2}) and the fact that $(\frac{2s}{p}+2t) \frac{p}{pq-1} = n-2s$, one gets
\[
C_1 \frac{\log(1+ |x|)}{|x|^{n-2s}} \leq m_{u}(r) \leq C_2 \frac{1}{|x|^{n-2s}}
\]
for all $x$ with $|x|=r$ large enough. But this contradicts the positivity of $u$.

It still remains the situation when $-n+2s \in (-1,0)$. In this case, the function $\Theta(r)$ is increasing near the origin and decreasing for $r$ large, with exactly one maximum point, say at $r_0 > 0$. Consider the function
\[
w(r) = \left\{
\begin{array}[c]{ll}%
\Theta(r_0) & {\rm if}\ 0 < r \leq r_0\\
\Theta(r) & {\rm if}\ r_0 < r
\end{array}
\right..
\]
Again, one defines the comparison function for $R_0 > 1$ and $R_0/2 > r_0$ as in Lemma \ref{lema4.2}
\[
\varphi(x) = m_{u}(R_0)\frac{w(r)-w(R_1)}{w(r_0)-w(R_1)}
\]
for $|x| < R_1$ and $\varphi(x)=0$ for $|x| \geq R_1$, where $R_1 > R_0$. It is clear that $\varphi(x) \leq u(x)$ for all $x$ with $|x| \leq R_0$ or $|x| \geq R_1$. In addition,
\[
(-\Delta)^{s} \varphi(x) \leq \frac{C}{|x|^{n}}
\]
for all $x$ with $R_0 < |x| < R_1$. Lastly, using Lemma \ref{lema6.1.} and the fact that $\Theta$ is increasing in $(0,r_0)$ and decreasing for $r \geq r_0$, the proof proceeds exactly as before and again we achieve the contradiction $u = 0$. This concludes the proof of sufficiency.\; \fim\\

\n {\bf Proof of the necessity of (\ref{4}).} Assume that the condition (\ref{4}) fails. In other words, we have

\begin{equation} \label{5}
\left(\frac{2s}{p}+2t\right)\frac{p}{pq-1} < n-2s\ \ {\rm and}\ \ \left(\frac{2t}{q}+2s\right)\frac{q}{pq-1} < n-2t\, .
\end{equation}

Consider the functions

\begin{equation}
u(x) = \frac{A}{(1+ |x|)^{2s k_{1}}} \ {\rm and}\ v(x) = \frac{B}{(1+ |x|)^{2t k_{2}}}\, ,
\end{equation}
where

\[
k_{1} = \frac{t+sp}{t(pq-1)}\ {\rm and}\ k_{2} = \frac{s+tq}{s(pq-1)}\, .
\]

The basic idea is to prove that $(u,v)$ is a positive radial super-solution of (\ref{sistema2}) with $G=\R^{n}$ for a suitable choice of positive constants $A$ and $B$.

Firstly, we assert that the inequalities

\begin{equation} \label{Ineq1}
\frac{1}{(1 - a + |ae_{1} + y|)^{2s k_{1}}} + \frac{1}{(1 - a + |ae_{1}-y|)^{2s k_{1}}}
\leq \frac{1}{|e_{1} + y|^{2s k_{1}}} + \frac{1}{|e_{1} - y|^{2s k_{1}}}
\end{equation}
\n and

\begin{equation} \label{Ineq2}
\frac{1}{(1-a+ |ae_{1} + y|)^{2t k_{2}}} + \frac{1}{(1-a+ |ae_{1} - y|)^{2t k_{2}}}
\leq \frac{1}{|e_{1} + y|^{2t k_{2}}} + \frac{1}{|e_{1} - y|^{2t k_{2}}}
\end{equation}

\n hold for all $a \in [0,1)$, $b \geq 0$ and $y \in \R$. In fact, consider the function $f(a,b,y)$ given by

\[
f(a,b,y) = (1-a+(a+b)^{2}+y^{2})^{1/2})^{-2 \alpha}+(1-a+(a-b)^{2}+y^{2})^{1/2})^{-2 \alpha}
\]

\[
- ((1+b)^{2}+y^{2})^{- \alpha}-((1-b)^{2}+y^{2})^{- \alpha}
\]
where $\alpha > 0$. One easily checks that

\begin{eqnarray*}
\frac{\partial f}{\partial a}(a,b,y) &=& \frac{-2 \alpha}{(1-a+(a+b)^{2}+y^{2})^{1/2})^{2 \alpha+1}} \left(-1+\frac{a+b}{((a+b)^{2}+y^{2})^{1/2}}\right)\\
&& +\frac{-2 \alpha}{(1-a+(a-b)^{2}+y^{2})^{1/2})^{2 \alpha +1}}\left(-1+\frac{a-b}{((a-b)^{2}+y^{2})^{1/2}}\right) \geq 0
\end{eqnarray*}
and $f(1,b,y) = 0$ for all $a \in [0,1)$, $b \geq 0$ and $y \in \R$. In particular, $f(a,b,y) \leq 0$ for all $a \in [0,1)$, $b \geq 0$ and $y \in \R$.

\n For $a=r/(1+r)$ and $x$ with $r = |x|$, we then have

\begin{eqnarray*}
&& \frac{1}{(1+ |x+y|)^{2 \alpha}} + \frac{1}{(1+ |x-y|)^{2 \alpha}} - \frac{2}{(1+ |x|)^{2 \alpha}}\\
& = & \frac{1}{(1+ |x|)^{2 \alpha}}\left\{\frac{1}{(1-a+ |ae_{1} + \overline{y}|)^{2 \alpha}} + \frac{1}{(1-a+ |ae_{1}-\overline{y}|)^{2 \alpha}} - 2\right\}\\
& \leq & \frac{1}{(1+ |x|)^{2 \alpha}}\left\{\frac{1}{|e_{1} + \overline{y}|^{2 \alpha}}+\frac{1}{|e_{1} -\overline{y}|^{2 \alpha}}-2\right\}\, ,
\end{eqnarray*}
where $\overline{y}= \frac{1}{1 + r}Py$, being $P$ an appropriate rotation matrix.

\n With the choice $\alpha = s k_1$ and $\alpha = t k_2$, we derive (\ref{Ineq1}) and (\ref{Ineq2}), respectively. Using these inequalities, we find
\begin{eqnarray*}
(-\Delta)^{s} u(x) &=&  -\frac{1}{2}\int\limits_{\R^{n}}\frac{A}{(1+ |x-y|)^{2s k_{1}} |y|^{n+2s}} + \frac{A}{(1+ |x+y|)^{2s k_{1}} |y|^{n+2s}} \\
&& - \frac{2A}{(1+ |x|)^{2s k_{1}} |y|^{n+2s}}\; dy\\
& \geq & - \frac{1}{2}\frac{A}{(1+ |x|)^{2s(k_{1} + 1)}} \int\limits_{\R^{n}}\frac{|e_{1} + y|^{-2s k_{1}}+ |e_{1}-y|^{-2s k_{1}}-2}{|y|^{n+2s}}\; dy \\
&=& \frac{ c_{1}A}{(1+ |x|)^{2s(k_{1} + 1)}}
\end{eqnarray*}
and
\begin{eqnarray*}
(-\Delta)^{t} v(x) &=& -\frac{1}{2}\int\limits_{\R^{n}}\frac{B}{(1+ |x-y|)^{2t k_{2}} |y|^{n+2t}} + \frac{B}{(1+ |x+y|)^{2t k_{2}} |y|^{n+2t}}  \\
&&- \frac{2B}{(1+ |x|)^{2t k_{2}} |y|^{n+2t}}\; dy\\
& \geq & - \frac{1}{2}\frac{B}{(1+ |x|)^{2t(k_{2} + 1)}} \int\limits_{\R^{n}} \frac{|e_{1}+y|^{-2t k_{2}} + |e_{1}-y|^{-2t k_{2}}-2}{|y|^{n+2t}}\; dy \\
&=& \frac{c_{2}B}{(1+ |x|)^{2t(k_{2} + 1)}}\, .
\end{eqnarray*}

\n Since $pq > 1$, there exist constants $k_1$ and $k_2$ such that $2s(k_{1} + 1) = 2tk_{2}p$ and $2t(k_{2}+1)=2sk_{1}q$. Thanks to (\ref{5}), it readily follows that $k_1$ and $k_2$ are positive, $2s k_{1} < n-2s$ and $2t k_{2} < n-2t$. These last two conditions guarantee the positivity of the above constants $c_{1}$ and $c_{2}$.

\n On the other hand, we have

\[
(-\Delta)^{s} u(x) - v^{p}(x) \geq \frac{c_{1}A}{(1+\vert x\vert)^{2s(k_{1}+1)}} - \frac{B^{p}}{(1+\vert x\vert)^{2tk_{2}p}} = \frac{c_{1}A - B^p}{(1+\vert x\vert)^{2s(k_{1}+1)}}
\]
and

\[
(-\Delta)^{t} v(x) - u^{q}(x) \geq \frac{c_{2}B}{(1+\vert x\vert)^{2t(k_{2}+1)}} - \frac{A^{q}}{(1+\vert x\vert)^{2sk_{1}q}} = \frac{c_{2}B - A^q}{(1+\vert x\vert)^{2t(k_{2}+1)}}
\]
for all $x \in \R^{n}$. Finally, the assumption $pq>1$ also allows us to choose $A=(c_{1}c_{2}^{p})^{\frac{1}{pq-1}}>0$ and $B=(c_{1}^{q}c_{2})^{\frac{1}{pq-1}}>0$ so that the right-hand side of the above inequalities are equal to zero. This concludes the proof of Theorem \ref{Teo2}.\; \fim\\

\section{Proof of Theorem \ref{Teo3}}

The first tool to be used in the proof of Theorem \ref{Teo3} is the following result whose proof is based on the method of moving plane.

\begin{propo}\label{prop4.1}
Let $(u,v)$ be a positive viscosity bounded solution of
\begin{equation}
\left\{\begin{array}{llll}
(-\Delta)^{s}u = v^p & {\rm in} \ \ \R^{n}_{+}\\
(-\Delta)^{t}v = u^q & {\rm in} \ \ \R^{n}_{+}\\
u=v=0 & {\rm in} \ \ \R^n\setminus\R^{n}_{+}
\end{array}\right.
\end{equation}
Assume $0 < s,t < 1$ and $p, q \geq 1$. Then, $u$ and $v$ are strictly increasing in $x_{n}$-direction.
\end{propo}

\n {\bf Proof of Proposition \ref{prop4.1}.} Let $\Sigma_{\mu}:=\{(\overline{x},x_{n})\in\R^{n}_{+}:\ 0 < x_{n} < \mu\}$ and $T_{\mu}:=\{(\overline{x},x_{n})\in\R^{n}_{+}:\ x_{n}=\mu\}$. For $x=(\overline{x},x_{n})\in\R^{n}$, we denote $u_{\mu}(x)=u(x_{\mu})$, $w_{\mu,u}(x)=u_{\mu}(x)-u(x)$, $v_{\mu}(x)=v(x_{\mu})$ and $w_{\mu,v}(x)=v_{\mu}(x)-v(x)$, where $\mu>0$ and $x_{\mu}=(\overline{x},2\mu -x_{n})$ for all $(\overline{x},x_{n})\in\R^{n-1}\times\R$. For any subset $A$ of $\R^{n},$ we write $A_{\mu}=\{x_{\mu}:\ x\in A\}$, the reflection of $A$ with respect to $T_{\mu}$.

We next divide the proof into two steps.

\n {\bf First step:} We here prove that if $\mu > 0$ is small enough, then $w_{\mu,u} > 0$ and $w_{\mu,v} > 0$ in $\sum_{\mu}$. For this purpose, we define

\[
\Sigma_{\mu,u}^{-}=\{x\in\Sigma_{\mu}:\ w_{\mu,u}(x)<0\}\text{ and }\Sigma_{\mu,v}^{-}=\{x\in\Sigma_{\mu}:\ w_{\mu,v}(x)<0\}.
\]
We first show that $\Sigma_{\mu,u}^{-}$ is empty if $\mu$ is small enough. Indeed, assume for a contradiction that $\Sigma_{\mu,u}^{-}$ is not empty and define

\begin{equation}  \label{4.2}
w_{\mu,u}^{1}(x)=\left\{\begin{array}{ccll}
w_{\mu,u}(x)& & {\rm if} \ \ x\in\Sigma_{\mu,u}^{-}\\
0&   & {\rm if} \ \ x \in \R^{n}\setminus\Sigma_{\mu,u}^{-}
\end{array}\right.
\end{equation}
and

\begin{equation}   \label{4.3}
w_{\mu,u}^{2}(x)=\left\{\begin{array}{ccll}
 0 & & {\rm if} \ \ x \in \Sigma_{\mu,u}^{-}\\
w_{\mu,u}(x) & & {\rm if} \ \ x \in \R^{n}\setminus\Sigma_{\mu,u}^{-}
\end{array}\right.
\end{equation}
It is clear that $w_{\mu,u}^{1}(x) = w_{\mu,u}(x) - w_{\mu,u}^{2}(x)$ for all $x\in\R^{n}$. For each $\mu>0$, we now assert that
\begin{equation} \label{3.4}
(-\Delta)^{s}w_{\mu,u}^{2}(x)\leq 0\text{ for all }x\in\Sigma_{\mu,u}^{-}\, .
\end{equation}
In fact, from the definition of $(-\Delta)^{s}$, we have

\begin{eqnarray*}
(-\Delta)^{s}w_{\mu,u}^{2}(x) &=& \int\limits_{\R^{n}}\frac{w_{\mu,u}^{2}(x)-w_{\mu,u}^{2}(y)}{ |x-y|^{n+2s}}\;dy = -\int\limits_{\R^{n}\setminus\Sigma_{\mu,u}^{-}}\frac{w_{\mu,u}^{2}(y)}{ | x-y|^{n+2s}}\; dy\\
&=& - \int\limits_{(\Sigma_{\mu}\setminus\Sigma_{\mu,u}^{-})\cup (\Sigma_{\mu}\setminus\Sigma_{\mu,u}^{-})_{\mu}}\frac{w_{\mu,u}(y)}{ |x-y|^{n+2s}}\; dy\\
&& - \int\limits_{(\R^{n}\setminus\R^{n}_{+})\cup(\R^{n}\setminus\R^{n}_{+})_{\mu}}\frac{w_{\mu,u}(y)}{ |x-y|^{n+2s}}\; dy - \int\limits_{(\Sigma_{\mu,u}^{-})_{\mu}}\frac{w_{\mu,u}(y)}{ | x-y|^{n+2s}}\; dy\\
&=& -A_{1}-A_{2}-A_{3}
\end{eqnarray*}
for all $x \in \Sigma_{\mu,u}^{-}$.

We next estimate separately each of these integrals.

\n Firstly, note that $w_{\mu,u}(y_{\mu}) = -w_{\mu,u}(y)$ for all $y\in\R^{n}$ and $w_{\mu,u}^{2}(y) \geq 0$ in $\Sigma_{\mu}\setminus\Sigma_{\mu,u}^{-}$. Then,

\begin{eqnarray*}
A_{1} &=& \int\limits_{(\Sigma_{\mu}\setminus\Sigma_{\mu,u}^{-})\cup (\Sigma_{\mu}\setminus\Sigma_{\mu,u}^{-})_{\mu}}\frac{w_{\mu,u}(y)}{ |x-y|^{n+2s}}\; dy\\
&=& \int\limits_{\Sigma_{\mu}\setminus\Sigma_{\mu,u}^{-}}\frac{w_{\mu,u}(y)}{|x-y|^{n+2s}}\; dy + \int\limits_{\Sigma_{\mu}\setminus\Sigma_{\mu,u}^{-}}\frac{w_{\mu,u}(y_{\mu})}{ |x-y_{\mu}|^{n+2s}}\; dy\\
&=& \int\limits_{\Sigma_{\mu}\setminus\Sigma_{\mu,u}^{-}}w_{\mu,u}(y)\left(\frac{1}{ |x-y|^{n+2s}}-\frac{1}{|x-y_{\mu}|^{n+2s}}\right)\; dy \geq 0\, ,
\end{eqnarray*}
since $ |x-y_{\mu}| > |x-y|$ for all $x \in \Sigma_{\mu,u}^{-}$ and $y \in \Sigma_{\mu}\setminus\Sigma_{\mu,u}^{-}$.

\n In order to discover the sign of $A_{2}$ we observe that $u=0$ in $\R^{n}\setminus\R^{n}_{+}$ and $u_{\mu}=0$ in $(\R^{n}\setminus\R^{n}_{+})_{\mu},$ so we have

\begin{eqnarray*}
A_{2}&=& \int\limits_{(\R^{n}\setminus\R^{n}_{+})\cup(\R^{n}\setminus\R^{n}_{+})_{\mu}}\frac{w_{\mu,u}(y)}{ |x-y|^{n+2s}}\; dy\\
&=& \int\limits_{\R^{n}\setminus\R^{n}_{+}}\frac{u_{\mu}(y)}{ |x-y|^{n+2s}}\; dy -\int\limits_{(\R^{n}\setminus\R^{n}_{+})_{\mu}}\frac{u(y)}{ |x-y|^{n+2s}}\; dy\\
&=& \int\limits_{\R^{n}\setminus\R^{n}_{+}}u_{\mu}(y)\left(\frac{1}{ |x-y|^{n+2s}}-\frac{1}{ |x-y_{\mu}|^{n+2s}}\right)\; dy \geq 0\, ,
\end{eqnarray*}
since $u_{\mu}\geq 0$ in $\R^{n} \setminus \R^{n}_{+}$ and $|x-y_{\mu}| > |x-y|$ for all $x \in \Sigma_{\mu,u}^{-}$ and $y \in \R^{n}\setminus\R^{n}_{+}$.

\n Finally, since $w_{\mu,u}<0$ in $\Sigma_{\mu,u}^{-}$, we have

\begin{eqnarray*}
A_{3} = \int\limits_{(\Sigma_{\mu,u}^{-})_{\mu}}\frac{w_{\mu,u}(y)}{ |x-y|^{n+2s}}\; dy = \int\limits_{\Sigma_{\mu,u}^{-}}\frac{w_{\mu,u}(y_{\mu})}{ |x-y_{\mu}|^{n+2s}}\; dy = -\int\limits_{\Sigma_{\mu,u}^{-}}\frac{w_{\mu,u}(y)}{ |x-y_{\mu}|^{n+2s}}\; dy \geq 0\, .
\end{eqnarray*}

\n Hence, the claim (\ref{3.4}) follows.

Using now (\ref{3.4}), for any $x \in \Sigma_{\mu,u}^{-}$, one has

\begin{eqnarray*}
(-\Delta)^{s}w_{\mu,u}^{1}(x)&=&(-\Delta)^{s}w_{\mu,u}(x)=(-\Delta)^{s}u_{\mu}(x)-(-\Delta)^{s}u(x)\\
&=& v^p_{\mu}(x) - v^p(x) = \frac{v^p_{\mu}(x) - v^p(x)}{v_{\mu}(x) - v(x)} w_{\mu,v}(x)\, .
\end{eqnarray*}

\n Define

\[
\varphi_{v}(x) = \frac{v^p_{\mu}(x) - v^p(x)}{v_{\mu}(x)-v(x)}
\]
for $x \in \Sigma_{\mu,u}^{-}$.

Since $p \geq 1$, we have $\varphi_{v} \in L^{\infty}(\Sigma_{\mu,u}^{-})$ and $\varphi_{v}w_{\mu,v}$ is continuous. In addition, since $w_{\mu,u}^{1} = 0$ in $\R^{n}\setminus\Sigma_{\mu,u}^{-}$, by Theorem 2.3 of \cite{xia}, one gets

\begin{equation} \label{est}
\Vert w_{\mu,u}^{1}\Vert_{L^{\infty}(\Sigma_{\mu,u}^{-})} \leq CR(\Sigma_{\mu,u}^{-})^{2s}\Vert \varphi_{v}w_{\mu,v}\Vert_{L^{\infty}(\Sigma_{\mu,u}^{-})}\, ,
\end{equation}
where $R(\Sigma_{\mu,u}^{-})$ is the smallest positive constant $R$ such that
\[
|B_{R}(x)\setminus\Sigma_{\mu,u}^{-}| \geq \frac{1}{2} |B_{R}(x)|
\]
for all $x \in \Sigma_{\mu,u}^{-}$. Besides, we have

\[
\varphi_{v}w_{\mu,v}(x) = v^p(x) - v^p_{\mu}(x) \leq 0 \ \ {\rm in}\ \ \Sigma_{\mu} \setminus \Sigma_{\mu,v}^{-}
\]
and
\[
\varphi_{v}w_{\mu,v}(x) = v^p(x) - v^p_{\mu}(x) > 0 \ \ {\rm in}\ \ \Sigma_{\mu,v}^{-}\, .
\]
Let $\Sigma_{\mu}^{-} = \Sigma_{\mu,u}^{-} \cap \Sigma_{\mu,v}^{-}$. Then, from (\ref{est}), one derives

\begin{eqnarray*}
\Vert w_{\mu,u}^{1}\Vert_{L^{\infty}(\Sigma_{\mu,u}^{-})}& \leq & CR(\Sigma_{\mu,u}^{-})^{2s}\Vert \varphi_{v}w_{\mu,v}\Vert_{L^{\infty}(\Sigma_{\mu}^{-})}\\
&\leq & CR(\Sigma_{\mu,u}^{-})^{2s}\Vert \varphi_{v}\Vert_{L^{\infty}(\Sigma_{\mu}^{-})}\Vert w_{\mu,v}\Vert_{L^{\infty}(\Sigma_{\mu}^{-})}\\
&\leq & CR(\Sigma_{\mu,u}^{-})^{2s}\Vert w_{\mu,v}\Vert_{L^{\infty}(\Sigma_{\mu}^{-})}\, ,
\end{eqnarray*}
where in the last inequality we use the condition $p \geq 1$.

Similar to (\ref{4.2}) and (\ref{4.3}), we define
\begin{equation}
w_{\mu,v}^{1}(x)=\left\{\begin{array}{ccll}
w_{\mu,v}(x) & & {\rm if} \ \ x\in\Sigma_{\mu,v}^{-}\\
0 & & {\rm if} \ \ x\in\R^{n}\setminus\Sigma_{\mu,v}^{-}
\end{array}\right.
\end{equation}
and

\begin{equation}
w_{\mu,v}^{2}(x)=\left\{\begin{array}{ccll}
 0 & & {\rm if} \ \ x\in\Sigma_{\mu,v}^{-}\\
w_{\mu,v}(x) & & {\rm if} \ \ x\in\R^{n}\setminus\Sigma_{\mu,v}^{-}
\end{array}\right.
\end{equation}
and argue in a completely analogous way with the aid of the assumption $q \geq 1$ to obtain

\[
\Vert w_{\mu,v}^{1}\Vert_{L^{\infty}(\Sigma_{\mu,v}^{-})} \leq CR(\Sigma_{\mu,v}^{-})^{2t}\Vert w_{\mu,u}\Vert_{L^{\infty}(\Sigma_{\mu}^{-})}\, .
\]
Thus,

\[
\Vert w_{\mu,u}^{1}\Vert_{L^{\infty}(\Sigma_{\mu,u}^{-})} \leq C^{2}R(\Sigma_{\mu,u}^{-})^{2s}R(\Sigma_{\mu,v}^{-})^{2t}\Vert w_{\mu,u}^{1}\Vert_{L^{\infty}(\Sigma_{\mu,u}^{-})}
\]
and

\[
\Vert w_{\mu,v}^{1}\Vert_{L^{\infty}(\Sigma_{\mu,v}^{-})} \leq C^{2}R(\Sigma_{\mu,u}^{-})^{2s}R(\Sigma_{\mu,v}^{-})^{2t}\Vert w_{\mu,v}^{1}\Vert_{L^{\infty}(\Sigma_{\mu,v}^{-})}\, .
\]
Now choosing $\mu$ small enough so that $C^{2}R(\Sigma_{\mu,u}^{-})^{2s}R(\Sigma_{\mu,v}^{-})^{2t} < 1$, we conclude that $\Vert w_{\mu,u}^{1}\Vert_{L^{\infty}(\Sigma_{\mu,u}^{-})} = 0$, so $\vert \Sigma_{\mu,u}^{-}\vert = 0$. Since $\Sigma_{\mu,u}^{-}$ is open, we deduce that $\Sigma_{\mu,u}^{-}$ is empty, which is a contradiction. Therefore, we get $w_{\mu,u}\geq 0$ in $\Sigma_{\mu}$ for $\mu > 0$ small enough. Similarly, one gets $w_{\mu,v}\geq 0$ in $\Sigma_{\mu}$ for $\mu > 0$ small enough too. Moreover, since the functions $u$ and $v$ are positive in $\R^{n}_{+}$ and $u=v=0$ in $\R^{n}\setminus\R^{n}_{+}$, it follows that $w_{\mu,u}$ and $w_{\mu,v}$ are positive in $\{x_{n}=0\}$ and then, by continuity, $w_{\mu,u} \neq 0$ and $w_{\mu,v}\neq 0$ in $\Sigma_{\mu}$.

In order to complete the proof of this step, we assert that if $w_{\mu,u} \geq 0$, $w_{\mu,v} \geq 0$, $w_{\mu,u} \neq 0$ and $w_{\mu,v} \neq 0$ in $\Sigma_{\mu}$ with $\mu>0$, then $w_{\mu,u} > 0$ and $w_{\mu,v}> 0$ in $\Sigma_{\mu}$. Indeed, we have

\[
(-\Delta)^{s}w_{\mu,u}(x) = v^p_{\mu}(x) - v^p(x) \geq 0\ \ {\rm in}\ \ \Sigma_{\mu}
\]
and

\[
(-\Delta)^{t}w_{\mu,v}(x) = u^q_{\mu}(x) - u^q(x) \geq 0\ \ {\rm in}\ \ \Sigma_{\mu}\, .
\]
Since $w_{\mu,u} \geq 0$, $w_{\mu,v} \geq 0$, $w_{\mu,u} \neq 0$ and $w_{\mu,v} \neq 0$ in $\Sigma_{\mu}$, by the Silvestre's strong maximum principle, the conclusion follows.\\

\n {\bf Second step:} Define
\[
\mu^* = sup\{\mu > 0:\ w_{\nu,u}>0,\ w_{\nu,v} > 0\ \ {\rm in}\ \ \Sigma_{\nu}\ \ {\rm for\ all}\ \ 0 < \nu < \mu\}\, .
\]
It is clear that $\mu^* > 0$ and $w_{\mu,u}>0$ and $w_{\mu,v}>0$ in $\Sigma_{\mu}$ for all $0 <\mu < \mu^*$, so that $u$ and $v$ are strictly increasing in $x_{n}$-direction. Indeed, for $0 < x_{n} < \overline{x}_{n} < \mu^*$, let $\mu = \frac{x_{n}+\overline{x}_{n}}{2}$. Since $w_{\mu,u} > 0$ and $w_{\mu,v} > 0$ in $\Sigma_{\mu}$, we have

\[
0<w_{\mu,u}(x',x_{n})=u_{\mu}(x',x_{n})-u(x',x_{n})=u(x',\overline{x}_{n})-u(x',x_{n})
\]
and

\[
0<w_{\mu,v}(x',x_{n})=v_{\mu}(x',x_{n})-v(x',x_{n})=v(x',\overline{x}_{n})-v(x',x_{n})\, ,
\]
so that $u(x',\overline{x}_{n}) > u(x',x_{n})$ and $v(x',\overline{x}_{n}) > v(x',x_{n})$, as claimed. Thus, the proposition is proved if we are able to show that $\mu^* = +\infty$.

Suppose for a contradiction that $\mu^*$ is finite. Now choose $\varepsilon_{0} > 0$ small enough such that the operators $(-\Delta)^{s}-\varphi_{v}$ and $(-\Delta)^{t}-\varphi_{u}$ satisfies the strong maximum principle in the open $\Sigma_{\mu^{\ast}+\varepsilon_{0}}\setminus\Sigma_{\mu^{\ast}-\varepsilon_{0}}$, see \cite{xia}. Here we use that $\varphi_{u}(x)=\frac{u_{\mu}^{q}(x)-u^{q}(x)}{u_{\mu}(x)-u(x)}$ and $\varphi_{v}(x)=\frac{v_{\mu}^{p}(x)-v^{p}(x)}{v_{\mu}(x)-v(x)}$) can be taken small in the $L^\infty$-norm, since $p, q > 1$. Therefore, $w_{\mu^* + \varepsilon_{0},u} > 0$ and $w_{\mu^* + \varepsilon_{0},v} > 0$ in $\Sigma_{\mu^* + \varepsilon_{0}}$, providing a contradiction.\; \fim\\

\begin{propo}\label{prop4.2}
Let $0 < s,t < 1$ and $p, q > 0$. If the system
\begin{equation}
\left\{\begin{array}{llll}
(-\Delta)^{s}u = v^p & {\rm in} \ \ \R^{n}_{+}\\
(-\Delta)^{t}v = u^q & {\rm in} \ \ \R^{n}_{+}\\
u=v=0 & {\rm in} \ \ \R^n\setminus\R^{n}_{+}
\end{array}\right.\label{sistema10}
\end{equation}
has a positive viscosity bounded solution, then the same system has a positive viscosity solution in $\R^{n-1}.$
\end{propo}

\n {\bf Proof of Proposition \ref{prop4.2}.} Let $(u,v)$ be a positive bounded solution of (\ref{sistema10}), that is there exists a constant $M$ such that $0 < u \leq M$ and $0 < v \leq M$ in $\R^n_+$. In the strip $\Sigma_{1}=\{x \in \R^{n}:\ 0 < x_{n} < 1\}$, we set

\[
u_k(x',x_{n}) = u(x',x_{n}+k)\ \ {\rm and}\ \ v_k(x',x_{n}) = v(x',x_{n}+k)\, .
\]
Note that $(u_k,v_k)$ solves the system (\ref{sistema10}) in $\Sigma_{1}$ for each integer $k \geq 1$. In addition, $0 < u_k \leq M$ and $0 < v_k \leq M$ in $\Sigma_{1}$. Thus,

\[
(-\Delta)^{s}u_k \leq M^{p}\ \ {\rm and}\ \ (-\Delta)^{s}u_k \geq 0\ \ {\rm in} \ \Sigma_{1}\, ,
\]

\[
(-\Delta)^{t}v_k \leq M^{q}\ \ {\rm and}\ \ (-\Delta)^{t}v_k \geq 0\ \ {\rm in} \ \Sigma_{1}\, .
\]
Then, by Theorem 2.6 of \cite{xia}, for any $\Omega' \subset \subset \Sigma_{1}$ and $0 < \beta < 1$, there exists a constant $C > 0$ such that $u_k, v_k \in C^{\beta}(\Omega')$ and

\[
\Vert u_k \Vert_{C^{\beta}(\Omega')} \leq C\left\{\Vert u_k\Vert_{L^{\infty}(\Sigma_{1})}+M^{p}\right\}
\]
and

\[
\Vert v_k \Vert_{C^{\beta}(\Omega')} \leq C\left\{\Vert v_k\Vert_{L^{\infty}(\Sigma_{1})}+M^{q}\right\}\, .
\]
So, the sequences $\{u_k\}$ and $\{v_k\}$ are bounded in $C^{\beta}(\Omega')$ and then, up to a subsequence, $\{u_k\}$ and $\{v_k\}$ converge uniformly on compact subset of $\Sigma_{1}$ to functions $\overline{u}$ and $\overline{v}$, respectively. By Theorem 2.7 of \cite{xia}, $(\overline{u},\overline{v})$ satisfies
\begin{equation}
\left\{\begin{array}{lc}
(-\Delta)^{s}\overline{u} = \overline{v}^{q} $ in $ \Sigma_{1} & \\
(-\Delta)^{t}\overline{v} = \overline{u}^{p} $ in $ \Sigma_{1}
\end{array}\right.   \label{sistema8}
\end{equation}

\n in the viscosity sense. The strict monotonicity provided in Proposition \ref{prop4.1} guarantees that $(\overline{u},\overline{v})$ is positive and independent of the $x_{n}$-variable.

\n On the other hand, the definition of $(-\Delta)^{s}$ gives
\begin{eqnarray*}
(-\Delta)^{s}\overline{u}(x)&=&\int\limits_{\R^{n-1}}\int\limits_{\R}\frac{\overline{u}(x')-\overline{u}(y')}{(|x'-y'|^{2}+(x_{n}-y_{n})^{2})^{\frac{n+2s}{2}}}\; dy_{n}\; dy'\\
&=&\int\limits_{\R^{n-1}}\int\limits_{\R}\frac{\overline{u}(x')-\overline{u}(x'-y')}{(|y'|^{2}+(y_{n})^{2})^{\frac{n+2s}{2}}}\; dy_{n}\; dy'\, .
\end{eqnarray*}
Let $y_{n}=|y'|\tan\theta$, where $\theta\in (-\frac{\pi}{2},\frac{\pi}{2})$, then
\begin{eqnarray*}
(-\Delta)^{s}\overline{u}(x) &=& \int\limits_{\R^{n-1}}\int\limits_{-\frac{\pi}{2}}^{\frac{\pi}{2}}\frac{\overline{u}(x')-\overline{u}(x'-y')}{|y'|^{n-1+2s}}(\cos\theta)^{n-2+2s}\; d\theta\; dy'\\
&=& \int\limits_{\R^{n-1}}\frac{\overline{u}(x') - \overline{u}(x'-y')}{|y'|^{n-1+2s}}\; dy' \int\limits_{-\frac{\pi}{2}}^{\frac{\pi}{2}}(\cos\theta)^{n-2+2s}\; d\theta
\end{eqnarray*}
and

\[
\int\limits_{-\frac{\pi}{2}}^{\frac{\pi}{2}}(\cos\theta)^{n-2+2s}\; d\theta = 2\int\limits_{0}^{\frac{\pi}{2}}(\cos\theta)^{n-2+2s}\; d\theta < + \infty\, ,
\]
since $n - 2 + 2s > 0$. This means that the $n$-dimension fractional Laplace operator is actually $(n-1)$-dimension, and we have
\begin{equation}
\left\{\begin{array}{lc}
(-\Delta)^{s}\overline{u} = \overline{v}^{q} $ in $ \R^{n-1} & \\
(-\Delta)^{t}\overline{v} = \overline{u}^{p} $ in $ \R^{n-1}
\end{array}\right.
\end{equation}																	\; \fim\\

Finally, Theorem \ref{Teo3} follows directly from Theorem \ref{Teo2} and Proposition \ref{prop4.2}.\; \fim\\

\section{Proof of Theorem \ref{teo1}}

The proof of the part of existence is an application of degree theory for compact operators in cones. This theory, essentially developed by Krasnoselskii, has often been used to show that certain operators admit fixed points. We are going to use an extension of Krasnoselskii results (se for instance \cite{sirakov}). The applicability of this theory relies on a priori bounds in $L^\infty$ of solutions of certain systems related to (\ref{1}) to be obtained through blow-up techniques by invoking Theorems \ref{Teo2} and \ref{Teo3}.

We begin by stating the above-mentioned abstract tool.

\begin{propo}\label{prop6.1}
Let $K$ be a closed cone with non-empty interior in a Banach space $X$ and let $T: K \rightarrow K$ and $H: [0,\infty) \times K \rightarrow K$ be continuous compact operators such that $T(0) = 0$ and $H(0,x) = T(x)$ for all $x\in K$. Assume there exist $\theta_0 > 0$ and $0 < r < R$ such that
    \begin{itemize}
    \item[(i)] $x \neq \theta T(x)$ for all $0 \leq \theta \leq 1$ and $x \in K$ such that $\Vert x\Vert =  r$,
    \item[(ii)] $H(\theta,x)\neq x$ for all $\theta \geq \theta_0$ and $x \in K$ with $\Vert x\Vert \leq R$,
    \item[(iii)] $H(\theta,x) \neq x$ for all $\theta \in [0,+\infty)$ and $x \in K$ with $\Vert x\Vert = R$.
    \end{itemize}
\n Then, $T$ has a fixed point $x_0 \in K$ such that $r \leq \Vert x_0\Vert \leq R$.
\end{propo}

Here $X$ denotes the Banach space $\{ (u,v) \in C(\R^n) \times C(\R^n):\ u,v=0\ {\rm in}\ \R^n\setminus\Omega\}$ endowed with the norm

\[
\Vert (u,v)\Vert := \max\{\Vert u\Vert_{L^{\infty}(\Omega)},\Vert v\Vert_{L^{\infty}(\Omega)}\}\, .
\]

\n and $K=\{u\in X:\ u,v\geq 0\ {\rm in}\ \Omega\}$. It is clear that solving (\ref{1}) is equivalent to finding a fixed point in $K$ of the operator $T:K \rightarrow K$ given by

\[
T(u,v)(x):=S(v^{p},u^{q})
\]
for $x \in \Omega$,where for any $(f,g) \in K$ we define $S(f,g)$ as the solution of the Dirichlet problem
\begin{equation}
\left\{\begin{array}{llll}
(-\Delta)^{s}u = f & {\rm in} \ \ \Omega\\
(-\Delta)^{t}v = g & {\rm in} \ \ \Omega\\
u=v=0 & {\rm in} \ \ \R^n\setminus\Omega
\end{array}\right.  \label{sistema6}
\end{equation}

\n Using that $\Omega$ is $C^2$ class, by Lemma 6.1 of \cite{sirakov}, the operator $S$ is well defined, linear, continuous and compact. Thus, one easily deduces that the operator $T$ is well defined, continuous and compact. In addition, we have $T(0,0) = 0$.

We also define $H: [0,\infty) \times K \rightarrow K$ as

\[
H(\theta,u,v) = S((v+\theta)^{p},(u+\theta)^{q})\, .
\]
Clearly, $H$ is well defined, continuous and compact too.

First we show that the condition (i) of Proposition \ref{prop6.1} is satisfied. This is the content of the following lemma:

\begin{lema} \label{lem1}
Assume that $0 < s, t < 1$ and $pq > 1$. Then, there exists a constant $r>0$ such that for any $\theta \in [0,1]$, the system
\begin{equation}
\left\{\begin{array}{llll}
(-\Delta)^{s}u = \theta v^p & {\rm in} \ \ \Omega\\
(-\Delta)^{t}v = \theta u^q & {\rm in} \ \ \Omega\\
u=v=0 & {\rm in} \ \ \R^n\setminus\Omega
\end{array}\right. \label{sistema7}
\end{equation}
has no classical solution $(u,v) \in K$ with $\Vert (u,v)\Vert = r$.
\end{lema}

\n {\bf Proof of Lemma \ref{lem1}.} We argue by contradiction. Let $\{(\theta_k,u_k,v_k)\}_{k\in\N}$ be a sequence of triples with $\theta_k \in [0,1]$ and $(u_k,v_k) \in K$ satisfying (\ref{sistema7}) such that $\Vert u_k\Vert_{L^{\infty}(\Omega)}, \Vert v_k\Vert_{L^{\infty}(\Omega)} \rightarrow 0$ as $k\longrightarrow +\infty$. Since $pq > 1$, we choose $\gamma$ such that

\[
\frac{1}{q} < \gamma < p
\]

\n and set $a_k = \Vert u_k\Vert_{L^{\infty}(\Omega)} + \Vert v_k\Vert^\gamma_{L^{\infty}(\Omega)}$. Define

\[
z_k=\frac{u_k}{a_k}\ \ {\rm and} \ \ w_k=\frac{v_k}{a_k^{1/\gamma}}\, .
\]
We then have

\[
(-\Delta)^{s}z_k = \frac{\theta_k}{a_k} v_k^p \ \ {\rm and} \ \ (-\Delta)^{t}w_k = \frac{\theta_k}{a_k^{1/\gamma}} u_k^q\, .
\]
Note that $\Vert z_k\Vert_{L^{\infty}(\Omega)} + \Vert w_k\Vert^\gamma_{L^{\infty}(\Omega)} = 1$,

\[
\left| \frac{\theta_k}{a_k} v_k^p \right| \leq \Vert v_k\Vert^{p - \gamma}_{L^{\infty}(\Omega)} \rightarrow 0\ \ {\rm and}\ \ \left| \frac{\theta_k}{a_k} u_k^q \right| \leq \Vert u_k\Vert^{q - 1/\gamma}_{L^{\infty}(\Omega)} \rightarrow 0
\]

\n uniformly for $x \in \Omega$. So, one easily deduces that $(z_k,w_k)$ converges uniformly to some couple $(z,w)$ satisfying $\Vert z\Vert_{L^{\infty}(\Omega)} + \Vert w\Vert^\gamma_{L^{\infty}(\Omega)} = 1$ and

\[
\left\{\begin{array}{llll}
(-\Delta)^{s}z = 0 & {\rm in} \ \ \Omega\\
(-\Delta)^{t}w = 0 & {\rm in} \ \ \Omega\\
z=w=0 & {\rm in} \ \ \R^n\setminus\Omega
\end{array}\right.
\]

\n But by uniqueness, we have $(z,w)=(0,0)$, providing a contradiction.\; \fim\\

The condition (ii) of Proposition \ref{prop6.1} follows from the following lemma:

\begin{lema} \label{lem2}
Assume that $0 < s, t < 1$, $p, q \geq 1$ and $pq > 1$. Then, there exists a constant $\theta_0>0$ such that for any $\theta \geq \theta_0$ the system
\begin{equation}
\left\{\begin{array}{llll}
(-\Delta)^{s}u = (v + \theta)^p & {\rm in} \ \ \Omega\\
(-\Delta)^{t}v = (u + \theta)^q & {\rm in} \ \ \Omega\\
u=v=0 & {\rm in} \ \ \R^n\setminus\Omega
\end{array}\right.  \label{sistema4}
\end{equation}
has no classical solution $(u,v)\in K$.
\end{lema}

\n {\bf Proof of Lemma \ref{lem2}.} Firstly, we define

\[
\lambda_1 := \inf \{ \int_\Omega |(-\Delta)^{s/2}u|^2 + |(-\Delta)^{t/2} v|^2\; dx:\ (u,v) \in H^s_0(\Omega) \times H^t_0(\Omega), \ \int_\Omega u^+ v^+\; dx = 1\}\, ,
\]
where $f^+ = \max\{f,0\}$. As usual, it follows that $\lambda_1$ is positive and attained for some couple $(\varphi, \psi) \in H^s_0(\Omega) \times H^t_0(\Omega)$. Also, by the weak maximum principle, $\varphi, \psi \geq 0$ in $\Omega$ and $\varphi, \psi \neq 0$ and, moreover, $(\varphi, \psi)$ satisfies

\[
\left\{\begin{array}{llll}
(-\Delta)^{s} \varphi = \lambda_1 \psi & {\rm in} \ \ \Omega\\
(-\Delta)^{t}\psi = \lambda_1 \varphi & {\rm in} \ \ \Omega\\
\varphi=\psi=0 & {\rm in} \ \ \R^n\setminus\Omega
\end{array}\right.
\]

On the other hand, by assumption, $p > 1$ or $q > 1$. If the first situation occurs, then for $A \geq \lambda_1^2$ there exists $\theta_0 > 0$ such that

\[
(y + \theta)^p \geq A (y + \theta) > A y\ \ {\rm and}\ \ (y + \theta)^p \geq (y + \theta) > y
\]
for all $y \geq 0$ and $\theta \geq \theta_0$.

Now let $\theta \geq \theta_0$ and $(u,v) \in K$ be a classical solution of (\ref{sistema4}). Then, by the Silvestre's strong maximum principle, we have $u, v > 0$ in $\Omega$ and

\[
\left\{\begin{array}{llll}
(-\Delta)^{s}u > A v & {\rm in} \ \ \Omega\\
(-\Delta)^{t}v > u & {\rm in} \ \ \Omega\\
u=v=0 & {\rm in} \ \ \R^n\setminus\Omega
\end{array}\right.
\]
Using the above equations satisfied by $(\varphi, \psi)$, one obtains

\[
\lambda_1 \int_\Omega u \psi\; dx > A \int_\Omega v \varphi\; dx\ \ {\rm and}\ \ \lambda_1 \int_\Omega v \varphi\; dx > \int_\Omega u \psi\; dx\; ,
\]
so that $A < \lambda_1^2$, providing a contradiction. \; \fim\\

Finally, the condition (iii) of Proposition \ref{prop6.1} is a consequence of the following lemma:

\begin{lema}\label{lem3}
Assume that $\Omega$ is of $C^2$ class, $0 < s, t < 1$, $n > 2s + 1$, $n > 2t + 1$, $p,q \geq 1$, $pq > 1$ and (\ref{4}) is satisfied. For each $\theta_0 > 0$ there exists a constant $C > 0$, depending only of $s,t,p,q$ and $\Omega$, such that for any classical solution $(u,v) \in K$ of the system (\ref{sistema4}) with $0 \leq \theta \leq \theta_{0}$, one has
\[
\Vert (u,v)\Vert \leq C\, .
\]
\end{lema}

\n {\bf Proof of Lemma \ref{lem3}.} Suppose for a contradiction that there exists a sequence $(u_k,v_k) \in K$ of solutions of (\ref{sistema4}) with $\theta = \theta_{k} \in [0,\theta_{0}]$ such that at least one of the sequence $(u_k)$ and $(v_k)$ tends to infinity in the $L^{\infty}$-norm.

Let $\beta_{1}=\left(\frac{2s}{p}+2t\right)\frac{p}{pq-1}$ and $\beta_{2}=\left(\frac{2t}{q}+2s\right)\frac{q}{pq-1}$. We set

\[
\lambda_k = \Vert u_k\Vert_{L^{\infty}(\Omega)}^{-\frac{1}{\beta_{1}}}\, ,
\]
if $\Vert u_k \Vert_{L^{\infty}(\Omega)}^{\beta_{2}} \geq \Vert v_k\Vert_{L^{\infty}(\Omega)}^{\beta_{1}}$, up to a subsequence, and $\lambda_k = \Vert v_k\Vert_{L^{\infty}(\Omega)}^{-\frac{1}{\beta_{2}}}$, otherwise. It suffices to assume the first of these two situations.

Note that $\lambda_k \rightarrow 0$ as $k \rightarrow + \infty$. Let $x_k \in \Omega$ be a maximum point of $u_k$. The functions

\[
z_k(x)=\lambda_k^{\beta_{1}}u_k(\lambda_k x + x_k)\ \ {\rm and}\ \ w_k(x)=\lambda_k^{\beta_{2}}v_k(\lambda_k x + x_k)
\]
are such that $z_k(0) = 1$ and $0 \leq z_k, w_k \leq 1$ in $\Omega_k:=\frac{1}{\lambda_k}(\Omega -x_k)$. Also, one checks that the functions $z_k$ and $w_k$ satisfy

\begin{equation}
\left\{\begin{array}{lc}
(-\Delta)^{s} z_k = \left(\lambda_k^{(2s+\beta_{1}-p\beta_{2})/p} w_k + \lambda_k^{(2s+\beta_{1})/p} \theta_k \right)^p = \left( w_k + \lambda_k^{(2s+\beta_{1})/p} \theta_k \right)^p\\
(-\Delta)^{t} w_k = \left(\lambda_k^{(2t+\beta_{2}-q\beta_{1})/q} z_k + \lambda_k^{(2t+\beta_{2})/q} \theta_k \right)^q = \left( z_k + \lambda_k^{(2t+\beta_{2})/q} \theta_k \right)^q
\end{array}\right. \label{36}
\end{equation}
in the open $\Omega_k$.

\n By compactness, module a subsequence, $(x_k)$ converges to some point $x_{0} \in \overline{\Omega}$. Let

\[
d_k = dist(x_k,\partial\Omega)\, .
\]
Two cases may occur as $k \rightarrow + \infty$:

\begin{itemize}
\item[(a)] $\frac{d_k}{\lambda_k} \rightarrow +\infty$, module a subsequence still denoted as before, or

\item[(b)] $\frac{d_k}{\lambda_k}$ is bounded.
\end{itemize}

\n If (a) occurs, then $\frac{1}{\lambda_k} B_{d_k}(0) \subset \Omega_k$ and $\frac{d_k}{\lambda_k} \rightarrow +\infty$ as $k \rightarrow + \infty$. So, $(\Omega_k)$ tends to $\R^{n}$ as $k \rightarrow +\infty$. We recall that $0 \leq z_k, w_k \leq 1$ in $\Omega_k$. Thus, the right-hand side of (\ref{36}) is bounded in $L^{\infty}(\Omega_k)$, so by compactness, we deduce that, up to a subsequence, $(z_k,w_k)$ converges to some function $(z,w)$ uniformly in compact sets of $\R^{n}$. By Theorem 2.7 of \cite{xia}, $(z,w)$ is a viscosity solution of (\ref{sistema2}) with $G=\R^{n}$. Note also that $z(0)=1$, since $z_k(0) = 1$ for all $k$, and hence $(z,w)\neq(0,0)$ and, by the Silvestre's strong maximum principle, $z, w > 0$ in $\R^n$. But this contradicts Theorem \ref{Teo2}.

Assume now that (b) occurs, that is $\frac{d_k}{\lambda_k}$ is bounded. In this case, up to a subsequence, we may assume that
\begin{equation}
\frac{d_k}{\lambda_k} \rightarrow a\in [0,\infty)\, .\label{6.11}
\end{equation}

\n Assume for a moment that $a > 0$. After a suitable rotation of $\R^{n}$ for each fixed $k$, one concludes that $(\Omega_k)$ converge to the half-space $\R^{n}_{+}=\{x \in \R^n: x_{n} > -a\}$. Again, we have $0 \leq z_k, w_k \leq 1$ in $\Omega_k$ and then, by compactness, $(z_k,w_k)$ converges, module a subsequence, to some function $(z,w)$ uniformly in compact sets of $\R^{n}_{+}$. As before, $(z,w)$ is a viscosity bounded solution of (\ref{sistema2}) with $G = \R^{n}_{+}$. Furthermore, using that $a > 0$ and $z_k(0)=1$ for all $k$, one gets $z(0)=1$, so that again $z,w > 0$ in $\Omega$ and this contradicts Theorem \ref{Teo3}.

The remainder of the proof consists in showing that $a>0$. We argue by contradiction and assume that $a=0$. The basic idea is to construct a barrier function $h_k$ on $\Omega_k$ for $z_k$. For this purpose, we define

\[
h_k(x)= (e^{-\frac{d_k}{\lambda_k}}-e^{x_{n}} )\sup\limits_{\Omega_k} \frac{( w_k + \lambda_k^{(2s+\beta_{1})/p} \theta_k )^p}{C_0}\, ,
\]
where $C_0$ is a positive constant such that

\begin{eqnarray*}
(-\Delta)^{s}e^{x_{n}}&=&-\int\limits_{\R^{n}}\frac{e^{(x_{n}+y_{n})}+e^{(x_{n}-y_{n})}-2e^{x_{n}}}{|y|^{n+2s}}\; dy\\
&=&-e^{x_{n}}\int\limits_{\R^{n}}\frac{e^{y_{n}}+e^{-y_{n}}-2}{|y|^{n+2s}}dy\leq -C_0 < 0
\end{eqnarray*}
for all $-\frac{d_k}{\lambda_k} < x_n < 0$. Thus, from (\ref{36}),
\[
(-\Delta)^{s}(h_k-z_k)\geq C_0\sup\limits_{\Omega_k} \left(\frac{( w_k + \lambda_k^{(2s+\beta_{1})/p} \theta_k )^p}{C_0}\right) - \frac{( w_k + \lambda_k^{(2s+\beta_{1})/p} \theta_k )^p}{C_0}\geq 0
\]
in $\Omega_k$ and $z_k \leq h_k$ in $\R^n\setminus\Omega_k$. Then, the weak maximum principle gives $z_k \leq h_k$ in $\Omega_k$. In addition, there exist $C_1 > 0$ and $\delta > 0$ such that

\[
|\nabla w_k(x)|\leq C_1
\]

\n for all $x \in \Omega_k \cap\{x\in\R^{n}:\ x_n + \frac{d_k}{\lambda_k} \leq \delta\}$. Since $x_k \in\Omega$, we have $0 \in \Omega_k \cap\{x\in\R^{n}:\ x_n + \frac{d_k}{\lambda_k} \leq \delta\}$ for $k$ large enough. Finally,

\[
1 = z_k(0) \leq h_k(0) \leq C_2 \left(e^{-\frac{d_k}{\lambda_k}}-1\right) \rightarrow 0
\]

\n as $k \rightarrow\infty$, providing a contradiction. \; \fim\\

\n Lastly, the conclusion of Theorem \ref{teo1} follows readily from Lemmas \ref{lem1}, \ref{lem2} and \ref{lem3} applied to Proposition \ref{prop6.1}.\\

\n {\bf Acknowledgments:} The first author was partially supported by CAPES and the second one was partially supported by CNPq and Fapemig. \\

 \end{document}